\documentclass[10pt]{amsart}
\let\newpf\proof \let\proof\relax 
\newenvironment{pf}{\newpf[\proofname]}{\qed\endtrivlist}

\def\even{{\rm even}}
\def\odd{{\rm odd}}

\def\be{\begin{equation}}
\def\ee{\end{equation}}

\def\ba{{\begin{align}}}
\def\ea{{\end{align}}}

\def\balpha{{\alpha_\theta}}

\def\0{{\mathbf 0}}

\def\cal{\mathcal}

\newtheorem{thm}{Theorem}[section]
\newtheorem{cor}[thm]{Corollary}

\newtheorem{lemma}[thm]{Lemma}

\newtheorem{prop}[thm]{Proposition}

\theoremstyle{remark}
\newtheorem{rem}{Remark}[section]

\numberwithin{equation}{section}

\def \bn {\hfill \\ \smallskip\noindent}

\theoremstyle{definition}

\def\proof{\bn {\bf Proof.} }

\def\note#1
{\marginpar
{\tiny $\leftarrow$
\par
\hfuzz=20pt \hbadness=9000 \hyphenpenalty=-100 \exhyphenpenalty=-100
\pretolerance=-1 \tolerance=9999 \doublehyphendemerits=-100000
\finalhyphendemerits=-100000 \baselineskip=6pt
#1}\hfuzz=1pt}

\newcommand{\LL}{{\cal L}}

\newcommand{\R}{{\mathbb R}}

\def\B0{{\bold{0}}}

\catcode`\@=12

\def\Empty{}
\newcommand\oplabel[1]{
  \def\OpArg{#1} \ifx \OpArg\Empty {} \else
  	\label{#1}
  \fi}
		
\newcommand{\comm}[1]{}
\newcommand{\comment}[1]{}

\begin{document}

\title{Convergence of an exact quantization scheme}

\author{Artur Avila}

\address{
Coll\`ege de France -- 3 Rue d'Ulm \\
75005 Paris -- France.
}
\email{avila@impa.br}

\comm{
\address{
IMPA -- Estr. D. Castorina 110 \\
22460-320 Rio de Janeiro -- Brazil.
}
\email{gugu@impa.br}
}

\thanks{Partially supported by Faperj and CNPq, Brazil.}   

\begin{abstract}

It has been shown by Voros \cite {V} that the spectrum
of the one-dimensional homogeneous anharmonic oscillator
(Schr\"odinger operator with potential $q^{2M}$, $M>1$)
is a fixed point of an explicit non-linear transformation.  We show that
this fixed point is globally and exponentially attractive in
spaces of properly normalized sequences.

\end{abstract}

\setcounter{tocdepth}{1}

\date\today

\maketitle

%\tableofcontents

\section{Introduction}

Let $0<\theta<\pi$ be a constant.  For $E,E'>0$, define
\be
\theta(E',E)=\tan^{-1} \frac {\sin \theta} {E'E^{-1}+\cos \theta}.
\ee

Let $X=(X_k)_{k=1}^\infty$, $Y=(Y_j)_{j=1}^\infty$ be sequences of positive
real numbers and define $\phi=(\phi_j)_{j=1}^\infty$ by
\be
\phi_j(X,Y)=\frac {1} {\pi} \sum_k \theta(X_k,Y_j).
\ee

Let $Q=(Q_i)_{i=1}^\infty$ be a constant vector, and consider the operator
$T \equiv T_{\theta,Q}$
given implicitly by $\phi(X,T(X))=Q$.
Of course $T(X)$ is only defined for certain sequences $X$.  We remark that
$T$ is dilatation equivariant ($T(\lambda X)=\lambda T(X)$ for $\lambda>0$)
and positive in the sense that if $0<X_k \leq X'_k$ for all $k>0$ and if
$T(X)=Y$ and $T(X')=Y'$ are defined then $Y_k \leq Y'_k$ for all $k>0$.

In this paper we will be interested in the description of
the dynamics of $T$ acting on certain spaces of normalized sequences,
under appropriate conditions on $Q$.

\subsection{Relation to exact anharmonic quantization}

We now describe the physical motivation of the problem (for futher details
and references, see \cite {V}, and for more recent related
work, see \cite {V1}).
Let us consider the one-dimensional anharmonic oscillator with even
homogeneous polynomial potential, that is, the Schr\"odinger operator
\be \label {oscillator}
(Hu)(q)=-\frac {d^2 u} {dq^2}+q^{2M} u(q), \quad M=2,3,\ldots,
\ee
acting on $L^2(\R)$.  This operator has a {\it discrete spectrum}
\be \label {spect}
0<E_0<E_1<\ldots,
\ee
where $\lim E_j=\infty$.

Let
\be \label {theta}
\theta=\frac {M-1} {M+1} \pi,
\ee
\be \label {alpha}
\balpha=\frac {\pi+\theta} {\pi}=\frac {2M} {M+1}.
\ee
It is known
that $E_k$ has polynomial growth, more precisely:

\begin{prop}[see \cite {V}, \S 2.1] \label {growthorder}

The spectrum (\ref {spect}) of the operator (\ref {oscillator}) satisfies
\be \label {nu}
\nu=\lim_{k \to \infty} k^{-\balpha} E_k,
\ee
where $\balpha$ is given by (\ref {alpha}) and $\nu$ is positive and finite.

\end{prop}

The semiclassical analysis provide much more information then what is
contained in the above proposition, for instance, $\nu$ can be explicitely
computed
\be
\nu=\left (2 \pi^{1/2} M \Gamma
\left (\frac {3} {2}+\frac {1} {2M} \right ) \Gamma \left (\frac {1} {2M}
\right )^{-1} \right )^\balpha,
\ee
and higher order terms for the asymptotic
development of $E_k$ are also available (though the resulting series does not
converge), let us only remark for motivation that
\be \label {alpha1}
E_k=\nu k^\balpha+O(k^{\balpha-1}).
\ee
It is convenient to split the spectrum according to parity
\be \label {P}
P^\even_i=E_{2i-2}, \quad P^\odd_i=E_{2i-1}, \quad i \geq 1.
\ee
It has been shown by Voros that $P^\even$ and
$P^\odd$ are fixed points of of operators
$T_{\theta,Q^\even}$ and $T_{\theta,Q^\odd}$ respectively, where
\be \label {Q}
Q^\even_k=k-\frac {3} {4}+\frac {M-1} {4(M+1)}, \quad
Q^\odd_k=k-\frac {1} {4}-\frac {M-1} {4(M+1)}.
\ee

\begin{prop} [see \cite {V}, \S 3.1] \label {fixedequation}

The even and odd parts of the spectrum of the operator (\ref {oscillator})
satisfy equations
\be
T_{\theta,Q^\even}(P^\even)=P^\even, \quad T_{\theta,Q^\odd}(P^\odd)=P^\odd,
\ee
where $\theta$, $Q^\even$ and $Q^\odd$ are as above.

\end{prop}

Due to dilatation equivariance of $T$,
the fixed point equation does not determine the spectrum completely. 
Numerical evidence was obtained (see \cite {V}, \S 7.3) that indicated
that it {\it does} determine the spectrum once one
normalizes appropriately at $k \to \infty$, and that the operator $T$
provides an exponentially convergent iterative scheme for determination
of the spectrum.

Our main theorem will confirm those hopes.  We show that there is only one
fixed point for $T$ subject to growth condition
$2^\balpha \nu k^\balpha+o(k^\balpha)$ (which thus
coincides with the spectrum $P$), and that
this is a globally attractive fixed point in the space of sequences with
such growth.  We also analyze the action of $T$ on sequences
whose growth is more accurately described in terms of polynomial error
terms, including the type
$2^\balpha \nu k^\balpha+O(k^{\balpha-\epsilon})$,
$0<\epsilon \leq 1$ (this is
natural in view of the asymptotic estimate (\ref {alpha1})), and we
show that the fixed point is indeed exponentially attractive among such
sequences.

\begin{thm} \label {spectrum}

Let $M>1$, and let $T$ denote the exact quantization operator related to the
even or odd spectrum $P$ of the operator (\ref {oscillator}).
If $X=(X_k)_{k=1}^\infty$ satisfies $X_k=2^\balpha
\nu e^{o(1)} k^\balpha$ then
$X^{(n)} \equiv T^n(X)$ converges pointwise to $P$, and indeed
\be
\lim_{n \to \infty} \sup_k k^{-\balpha}|X^{(n)}_k-P_k|=0.
\ee
If moreover $X_k=P_k+O(k^{\balpha-\epsilon})$
with $0<\epsilon<2$ then
\be
\sup_k k^{-\balpha+\epsilon}|X^{(n)}_k-P_k| \leq C \lambda^n,
\ee
where $C=C(X,\epsilon)>0$ and $\lambda=\lambda(\epsilon)<1$.

\end{thm}

The operator $T$ actually comes about as a (non-linear) quantization of a
semiclassical Bohr-Sommerfeld linear operator.  The main steps of our
analysis involves showing that $T$ behaves as a perturbation of the linear
operator.  The asymptotic limit $k \to \infty$ is given to certain accuracy
by the semiclassical linear operator, which can be shown to have the
required properties.  We use two obvious features of $T$ to show that the
quantization does not destroy those properties.  The first one is
positivity, and the second one is equivariance by dilatation.  Those
properties are present both at the infinitesimal analysis (they are used in
perturbative estimates of the operator norm of the derivative
$DT$) as in the
global analysis (where they are used in a key precompactness argument).

\section{Proof of Theorem \ref {spectrum}}

\subsection{Setting and notations}

We will actually prove a slightly more general result,
Theorem \ref {actual result},
about the operators $T_{\theta,Q}$.  This result implies Theorem \ref
{spectrum} immediately, using Propositions \ref {growthorder} and \ref
{fixedequation}.
%some of the above mentioned results
%from \cite {V}, namely the validity of the fixed point representation for
%the even and odd parts of the spectrum and that there exists $0<\nu<\infty$
%and $\alpha>0$ for which the asymptotic estimate
%(\ref {nu}) is valid.
The remaining analysis is completely self-contained.

We will need to make no restriction on $0<\theta<\pi$.
We will make two assumptions on the sequence $Q_k$:
\be \label {condition}
Q_k=k+O(1)
\ee
\be \label {condition'}
Q_k>\left (k-\frac {1} {2} \right ) \frac {\theta} {\pi}.
\ee
The first condition comes from the physical problem, and can be
relaxed to $Q_k=e^{o(1)} k$ without any changes in our analysis.
Notice that for any $X$,
\be
\sum_{j=1}^k \phi_j(X,X)>\frac {\theta k^2} {2 \pi} \quad k \geq 1,
\ee
in particular, if $\sum_{j=1}^k Q_j \leq (2 \pi)^{-1}\theta k^2$ for some
$k$ then there is no fixed point for $T$,
so some condition (possibly weaker)
in the line of our second condition is necessary for our results to hold.

It will be convenient to work in logarithmic coordinates for computations. 
All variables in capital letters will denote positive real numbers (or
vectors of positive real numbers).  The corresponding non-capital letters
will be reserved for their logarithms.

\subsection{Some spaces of sequences}

Let $u(\epsilon)$ be the space of $v=(v_i)_{i=1}^\infty$ of the form
\be
v_k=O(k^{-\epsilon}),
\ee
with the norm
\be
\|v\|_\epsilon=\sup k^\epsilon |v_k|.
\ee
Let $u^0(\epsilon)$ be the subspace of $u(\epsilon)$ consisting of
$v$ of the form
\be
v_k=o(k^{-\epsilon}).
\ee

Given some vector $x$, we define affine spaces
\be
u(x,\epsilon)=x+u(\epsilon)
\ee
and
\be
u^0(x,\epsilon)=x+u^0(\epsilon).
\ee
We will use the special notation
\be
u(\alpha,\epsilon)=u((\alpha \ln k)_{k=1}^\infty,\epsilon), \quad \alpha>0,
\epsilon \geq 0,
\ee
\be
u^0(\alpha,\epsilon)=u^0((\alpha \ln k)_{k=1}^\infty,\epsilon), \quad
\alpha>0, \epsilon \geq 0.
\ee
Notice that if $x \in u(\alpha,\epsilon)$ then
$u(x,\epsilon')=u(\alpha,\epsilon')$ provided $\epsilon' \leq \epsilon$.

The several affine spaces $u$ parametrize by exponentiation spaces $U$, for
instance
\be
U(\alpha,\epsilon)=\{(X_k)_{k=1}^\infty, X_k>0,
X_k=k^\alpha+O(k^{\alpha-\epsilon})\},
\ee
\be
U^0(\alpha,0)=\{(X_k)_{k=1}^\infty, X_k>0,
X_k=k^\alpha+o(k^\alpha)\}.
\ee

We can now state our main result:

\begin{thm} \label {actual result}

There exists a unique $\balpha>0$ for which there
exists a fixed point $X \in U(\balpha,0)$ for $T$.
Moreover,
\begin{enumerate}
\item The space $U^0(\balpha,0)$ is invariant for $T$,
\item There exists a fixed point $P \in U(\balpha,1)$,
\item $P$ is a global attractor in $U^0(\balpha,0)$, that is, for any $X \in
U^0(\balpha,0)$,
\be
\lim_{n \to \infty} \|T^n(x)-p\|_0=0,
\ee
\item The spaces $U(P,\epsilon)$ are invariant for $0 \leq
\epsilon<\balpha+1$,
\item $P$ is a global exponential attractor in $U(P,\epsilon)$,
$0<\epsilon<2$, that is, for any $X \in U(P,\epsilon)$,
\be
\|T^n(x)-p\|_\epsilon \leq C \lambda^n,
\ee
where $C=C(\epsilon,\|x-p\|_\epsilon)>0$ and $\lambda=\lambda(\epsilon)<1$.
\end{enumerate}

\end{thm}

The proof of this result will take the remaining of this section.

\subsection{Lipschitz continuity in $U(X,0)$}

Let us write $X \leq X'$ if $X_k \leq X_k'$ for all $k$.
Then $X \leq X'$ and $Y \geq Y'$ implies $\phi(X,Y) \geq \phi(X',Y')$, which
implies the positivity of $T$ we stated before:
$X \leq X'$ implies $T(X) \leq T(X')$.
In particular, $T(X) \leq X$ if $\phi(X,X) \geq Q$ and $T(X) \geq X$ if
$\phi(X,X) \leq Q$.

This also gives us a way to show that $T$ is defined at some $X$:
if $\phi(X,\underline Y) \leq Q \leq \phi(X,\overline Y)$ then
$T(X)=Y$ is defined and $\underline Y \leq Y \leq \overline Y$.

\begin{lemma} \label {lipschtz}

Assume that $T(X)=Y$ is defined.  Then $T$ is defined on $U(X,0)$ and
$T(U(X,0))=U(Y,0)$.  Moreover, $T:U(X,0) \to U(Y,0)$ is $1$-Lipschitz.

\end{lemma}

\begin{pf}

If $C^{-1}X \leq X' \leq CX$ then $\phi(X',C^{-1}Y) \leq
\phi(C^{-1}X,C^{-1}Y)=Q=\phi(CX,CY) \leq \phi(X',CY)$.
\end{pf}

\subsection{The derivative}

Let
\be
P(E,E')=\frac {E E'} {E^2+2\cos \theta E E'+E'^2}.
\ee
Notice that
\be
\frac {d\phi_j} {dx_k}(X,Y)=\frac {-\sin \theta} {\pi} P(X_k,Y_j),
\ee
\be
\frac {d\phi_j} {dy_j}(X,Y)=\sum_k \frac {\sin \theta} {\pi} P(X_k,Y_j),
\ee
and of course
\be
\frac {d\phi_j} {dy_k}(X,Y)=0, \quad j \neq k.
\ee

\comm{
Let us consider a curve $T(X(t))=Y(t)$.  It follows
\be
\frac {dy_i} {dx_j}(t)=\frac {P(X_j(t),Y_i(t))} {\sum_j
P(X_j(t),Y_i(t))}.
\ee 
}

We can now use write a nice formal expression for the derivative of
$T$ with respect to logarithmic coordinates.
If $T(x)=y$ is defined, let $DT(x)=(D_{ij}T(x))_{i,j \geq 1}$
be the infinite matrix
\be
D_{ij}T(x)=\frac {P(X_j,Y_i)} {\sum_k P(X_k,Y_i)}.
\ee
This matrix is stochastic and positive, that is all entries are positive
numbers and the sum of the entries in each row is $1$.  In particular,
the operator norm of $DT$ acting on bounded sequences is equal to $1$.

\begin{lemma} \label {lipderiv}

Let $\LL(u(0),u(0))$ be the space of bounded
linear transformations on $u(0)$ with the operator norm.
If $T(X)=Y$ is defined then $DT:u(x,0) \to \LL(u(0),u(0))$ is $4$-Lipschitz.

\end{lemma}

\begin{pf}

It follows immediately from the fact that $T$ is $1$-Lipschitz in $u(x,0)$
that if $\|x'-x''\|_0 \leq C$ then for all $i,j>0$,
\be
e^{-4C} \leq \frac {D_{ij}T(x')} {D_{ij}T(x'')} \leq e^{4C},
\ee
which easily implies the result.
\end{pf}

Notice that the previous proof implies that
\be
\|T(x+v)-T(x)-DT(x)v\|_0 \leq 4\|v\|_0^2,
\ee
so $DT$ is the actual derivative of $T:u(x,0) \to u(T(x),0)$.

\subsection{Weak contraction of $DT$ in $u^0(x,0)$}

\begin{lemma} \label {ineq}

Let $T(X)=Y$ be defined.  If $0 \neq v \in u^0(0)$ then
$\|DT(x)v\|_0<\|v\|_0$.

\end{lemma}

\begin{pf}

This is automatic since $DT$ is a stationary positive matrix.
\end{pf}

\begin{cor}

If $T(X)=Y$ is defined and $0 \neq v \in u^0(0)$ then
$\|T(x+v)-T(x)\|_0<\|v\|_0$.

\end{cor}

\begin{pf}

Integrate the previous estimate.
\end{pf}

\begin{cor}

There exists at most one fixed point in each $U^0(x,0)$.

\end{cor}

\subsection{The drift}

Let us define the drift
\be
D_\alpha=
\frac {1} {\pi} \int_0^\infty \tan^{-1} \frac {\sin \theta} {s^\alpha+\cos
\theta} ds.
\ee

\begin{lemma} \label {drif}

The operator $T$ is defined in $U(\alpha,0)$ if and only if $\alpha>1$, and
in this case the spaces $U(\alpha,0)$ are invariant.
Moreover, if $X \in U(\alpha,0)$,
then letting $T^n(X)=X^{(n)}$ we have
\be
\liminf_{k \to \infty} x_k-\alpha \ln k \leq \liminf_{k \to \infty}
x^{(n)}_k-\alpha \ln k+n\alpha \ln D_\alpha,
\ee
\be
\limsup_{k \to \infty} x^{(n)}_k-\alpha \ln k+n\alpha \ln D_\alpha \leq
\limsup_{k \to \infty} x_k-\alpha \ln k.
\ee

\end{lemma}

\begin{pf}

Let $T(X)=Y$, with $x_k \leq \alpha \ln k+C+o(1)$.  Then a simple
computation gives
\be
\phi_j(X,Y) \geq e^{-C\alpha^{-1}+o(1)} D_\alpha Y_j^{1/\alpha},
\ee
and since $\phi_j(X,Y)=j+O(1)$, we have
$Y_j \leq e^{C+o(1)} D_\alpha^{-\alpha} j^\alpha$.
Analogously, if $x_k \geq \alpha \ln k+C+o(1)$ then
\be
\phi_j(X,Y) \leq e^{-C\alpha^{-1}+o(1)} D_\alpha Y_j^{1/\alpha},
\ee
and since $\phi_j(X,Y)=j+O(1)$, we have
$Y_j \geq e^{C+o(1)} D_\alpha^{-\alpha} j^\alpha$.
\end{pf}

In particular, if $X \in U(\alpha)$ then the iterates of $X$ drift
(pointwise) towards either $0$ or $\infty$ unless $D_\alpha=1$.
Notice that
\be
\frac {d} {d\alpha} D_\alpha=\frac {-1} {\pi} \int_0^\infty
\frac {\sin \theta \ln s} {s^\alpha+2 \cos \theta+s^{-\alpha}} ds=\frac {-1}
{\pi} \int_1^\infty \frac {\sin \theta \ln s (1-s^{-2})} {s^\alpha+2 \cos
\theta+s^{-\alpha}} ds<0,
\ee
\be
\lim_{\alpha \to 1} D_\alpha=\infty,
\quad \lim_{\alpha \to \infty} D_\alpha=\frac {\theta} {\pi},
\ee
thus there exists a unique $\balpha>1$ such that $D_\balpha=1$.
From now on, $\balpha$ will denote this precise value.

\begin{rem} \label {r1}

One can actually compute explicitly
\be
D_\alpha=\frac {\sin \left (\frac {\theta} {\alpha} \right )} {\sin \left
(\frac {\pi} {\alpha} \right )},
\ee
so $D_\balpha=1$ implies $\balpha=1+\frac {\theta} {\pi}$.

\end{rem}

\begin{cor}

The space $U^0(\balpha,0)$ is invariant.

\end{cor}

\subsection{Construction of invariant sets}

Let $U$ be one of the spaces defined.  We say that
$K \subset U$ is uniformly bounded in $U$
if there exists $\underline X \leq \overline X$ in $U$ such that for all $Y
\in K$, $\underline X \leq Y \leq \overline X$.  Notice that the notion of
uniformly bounded in $U^0(\balpha,0)$ coincides with precompactness, while
the notion of uniformly bounded in $U(\balpha,0)$ coincides with
``bounded diameter''.

\begin{lemma}

In this setting

\begin{enumerate}

\item There exists $\overline X \in U(\balpha,1)$ with
$T(\overline X) \leq \overline X$, and $\overline X$
can be chosen arbitrarily big,

\item There exists $\underline X \in U(\balpha,1)$ with
$T(\underline X) \geq \underline X$ and $\underline X$ can be chosen
arbitrarily small.

\end{enumerate}

\end{lemma}

\begin{pf}

(Here, more precisely in the proof of (2),
is the only time we will use the condition (\ref {condition'}) on $Q_k$.)

(1)\, Let $Q_k<k+K$.  The required $\overline X$ is given by
$\overline X_k=(k+A)^\balpha$ for
all $A$ sufficiently big.  To see this, we must estimate,
for $A$ sufficiently big
\be
\phi_j(\overline X,\overline X)>j+K
\ee
for all $j$.

One can approximate
\be
\phi_j(\overline X,\overline X)=
\frac {1} {\pi} \int_A^\infty \tan^{-1} \frac {\sin \theta}
{s^\balpha (j+A)^{-\balpha}+\cos \theta} ds+O(1)
\ee
where the $O(1)$ term does not depend on $A$.
Of course
\be
\frac {1} {\pi} \int_A^\infty \tan^{-1} \frac {\sin \theta}
{s^\balpha (j+A)^{-\balpha}+\cos \theta} ds=
(j+A) \frac {1} {\pi} \int_{\frac {A} {j+A}}^\infty \tan^{-1}
\frac {\sin \theta} {s^\balpha+\cos \theta} ds.
\ee
This last term can be rewritten (using the condition on $\balpha$) as
\be
(j+A)-(j+A) \frac {1} {\pi} \int_0^{\frac {A} {j+A}} \tan^{-1}
\frac {\sin \theta} {s^\balpha+\cos \theta} ds.
\ee

Let us show the inequality (which trivially implies the required bound)
\be
j+\left (1-\frac {\theta} {\pi} \right ) A \leq
(j+A)-(j+A) \frac {1} {\pi} \int_0^{\frac {A} {j+A}} \tan^{-1}
\frac {\sin \theta} {s^\balpha+\cos \theta} ds \leq j+A,
\ee
or equivalently, with $B=(j+A)/A$,
\be
1-\frac {\theta} {\pi} \leq
1-B \frac {1} {\pi} \int_0^{B^{-1}} \tan^{-1} \frac {\sin \theta}
{s^\balpha+\cos \theta} ds \leq 1.
\ee
The right inequality being trivial, we estimate the left one
\be
\frac {\theta} {\pi} \geq
B \frac {1} {\pi} \int_0^{B^{-1}} \tan^{-1} \frac {\sin \theta}
{s^\balpha+\cos \theta} ds,
\ee
which is obvious since the integrand is a decreasing function of $s$ which
tends to $\theta$ when $s$ tends to $0$.

(2)\, Let $Q_k>k-K$.  The required $\underline X$ is given by
\begin{align}
\underline X_k&=N^{k-N^2}, &&k<N^2,\\
\underline X_k&=(k-N^2+N)^\balpha, && k \geq N^2,
\end{align}
for $N$ sufficiently big.
We can estimate, as before, for $j \geq N^2$
\be
\phi_j(\underline X,\underline X)=
\frac {1} {\pi} \sum_{k=1}^\infty \theta(\underline X_k,\underline X_j)=
\frac {1} {\pi} \sum_{k=1}^{N^2-1} \theta(\underline X_k,\underline X_j)+
\frac {1} {\pi} \sum_{k=N^2}^\infty \theta(\underline X_k,\underline
X_j)<j-K,
\ee
since
\be
\frac {1} {\pi} \sum_{k=1}^{N^2-1} \theta(\underline X_k,\underline X_j)
\leq (N^2-1) \frac {\theta} {\pi},
\ee
\be
\frac {1} {\pi} \sum_{k=N^2}^\infty \theta(\underline X_k,\underline X_j)
\leq j-N^2+N+O(1)
\ee
(the $O(1)$ independent of $N$ and $j$).
For $1 \leq j<N^2$ we estimate
\be
\phi_j(\underline X,\underline X)=\frac {1} {\pi} \sum_{k=1}^\infty
\theta(\underline X_k,\underline X_j)=
\frac {1} {\pi} \sum_{k=1}^j \theta(\underline X_k,\underline X_j)+
\frac {1} {\pi} \sum_{k=j+1}^\infty
\theta(\underline X_k,\underline X_j)<j-K
\ee
(implying the result), since
\be
\frac {1} {\pi} \sum_{k=1}^j \theta(\underline X_k,\underline X_j) \leq
j \frac {\theta} {\pi}-\frac {\theta} {2 \pi}
\ee
and
\be
\frac {1} {\pi} \sum_{k=j+1}^\infty \theta(\underline X_k,\underline
X_j)=o(1).
\ee
(the $o(1)$ in terms of $N$ and independent of $j$).
\end{pf}

\begin{cor} \label {cor1}

There exists a fixed point $P \in U(\balpha,1)$.  For any initial
condition $Y \in U(\balpha,1)$, $T^n(Y)$ converges to $P$ in
$U^0(\balpha,0)$.

\end{cor}

\begin{pf}

The previous lemma gives us
$\overline X, \underline X \in U(\balpha,1)$
with $\underline X \leq Y \leq \overline X$ and
with $\underline X \leq T(\underline X) \leq
T(\overline X) \leq \overline X$.  It follows that $T^n(\overline X)$
decreases pointwise to some vector
$\underline X \leq P \leq \overline X$.
This vector is obviously a fixed point of $T$.  This proves existence of the
fixed point.  Analogously,
$T^n(\underline X)$ increases pointwise to some fixed point, which must be
the same by uniqueness.  In particular, $T^n(Y)$ converges to $P$.
\end{pf}

\begin{lemma} \label {lemma2}

Let $Y \in U^0(\balpha,0)$.  Then $T^n(Y) \to P$ in the
$U^0(\balpha,0)$ metric.

\end{lemma}

\begin{pf}

We must show that for any $Y$, for any $\epsilon>0$, there exists $n_0$ such
that for $n>n_0$, $\|T^n(y)-p\|_0<\epsilon$.  Using the previous
construction, we obtain vectors $\underline X, \overline X
\in U(\balpha,1)$ with $(1-\epsilon/3)\underline X \leq Y \leq
(1+\epsilon/3) \overline X$.  Let $n_0$ be such that
\be
T^{n_0}(\overline X) \leq (1+\epsilon/3) P,
\ee
\be
T^{n_0}(\underline X) \geq (1-\epsilon/3) P.
\ee
It follows that for any $n>n_0$,
\be
(1-\epsilon/3)^2 P \leq T^n(Y) \leq (1+\epsilon/3)^2 P
\ee
which gives the desired estimate.
\end{pf}

\begin{cor} \label {precompact}

Let $K \subset U^0(\balpha,0)$ be uniformly bounded.  Then $\cup_{n=0}^\infty
T^n(K)$ is uniformly bounded as well.

\end{cor}

\begin{pf}

Let $\underline X \leq Y \leq \overline X$ for all $Y \in K$.  Then
$T^n(\underline X), T^n(\overline X) \to P$ implies that $\{T^n(\underline
X), T^n(\overline X)\}_{n \geq 0}$ is precompact, so uniformly bounded, by
say $\underline X'$, $\overline X'$.  By positivity of $T$, for any $Y \in
K$ one has $\underline X' \leq T^n(Y) \leq \overline X'$ as well.
\end{pf}

\subsection{Strong contraction of $DT$}

Let
\be
S_\epsilon=\int_0^\infty
\frac {s^{-\epsilon}} {s^\balpha+2 \cos \theta+s^{-\balpha}} ds.
\ee
Notice that
\be
S_\epsilon=\int_1^\infty \frac {s^{-\epsilon}+s^{\epsilon-2}}
{s^\balpha+2 \cos \theta+s^{-\balpha}} ds,
\ee
so that if $|\epsilon-1| \geq \balpha$ then $S_\epsilon=\infty$ and if
$|\epsilon-1|<\balpha$ then
$S_\epsilon=S_{2-\epsilon}$ is a strictly increasing function of
$|\epsilon-1|$.

\begin{rem} \label {r2}

It is possible to compute explicitly
\be
S_\epsilon=\frac {\pi} {\balpha \sin \theta}
\frac {\sin \left ((1-\epsilon) \theta \balpha^{-1} \right )}
{\sin \left ((1-\epsilon) \pi
\balpha^{-1} \right )}, \quad 0<|\epsilon-1|<\alpha_0,
\ee
while $S_1=\lim_{\epsilon \to 1} S_\epsilon=\frac {\theta}
{\alpha_\theta \sin \theta}$.  In particular, using that
$\alpha_\theta=1+\frac {\theta} {\pi}$ (Remark~\ref {r1}) we get
$S_1=\frac {\pi} {\alpha_\theta \sin \theta}$.

\end{rem}

\begin{lemma} \label {strong}

Let $K$ be a uniformly bounded set in $U^0(\balpha,0)$.
If $|\epsilon-1| \geq \balpha$ then for every
$X \in K$ we have that $DT(X)$ is not a bounded operator in
$u(\epsilon)$.  If $|\epsilon-1|<\balpha$ then there exists a norm $\| \cdot
\|_c$ in $u(\epsilon)$ (equivalent to $\| \cdot \|_\epsilon$)
and a constant $C_\epsilon$ such that $\|DT(X) v\|_c
\leq C_\epsilon \|v\|_c$ for $v \in u(\epsilon)$. 
Moreover, $C_\epsilon<1$ for $|\epsilon-1|<2$.

\end{lemma}

\begin{pf}

For $X \in K$, $X$ and $T(X)=Y$ satisfy uniformly
$x_k,y_k=\balpha \ln k+o(1)$.
Let $v_k=k^{-\epsilon}$ and let $w=DT(X) v$.  We have
\be
w_j=\left (\sum_k \frac {X_k Y_j} {X_k^2+2 \cos \theta X_k Y_j+Y_j^2}
v_k \right ) \left (
\sum_k \frac {X_k Y_j} {X_k^2+2 \cos \theta X_k Y_j+Y_j^2} \right )^{-1}.
\ee
which can be estimated as
\be
w_j=\left (\sum_k \frac {e^{o_{\min \{j,k\}}(1)}
k^\balpha j^\balpha} {k^{2\balpha}+
2 \cos \theta k^\balpha j^\balpha+j^{2\balpha}} k^{-\epsilon} \right )
\left (\sum_k \frac {e^{o_{\min\{j,k\}}(1)}
k^\balpha j^\balpha} {k^{2\balpha}+2 \cos \theta k^\balpha
j^\balpha+j^{2\balpha}} \right )^{-1}.
\ee
We easily estimate
\be
\sum_k e^{o_{\min\{j,k\}}(1)}
\frac {k^\balpha j^\balpha} {k^{2\balpha}+2 \cos \theta k^\balpha
j^\balpha+j^{2\balpha}}=
%e^{o_k(1)} \int_0^\infty \frac {j^\balpha t^\balpha} {t^{2\balpha}+2 \cos
%\theta t^\balpha j^\balpha+j^{2\balpha}} dt=
e^{o_j(1)} j S_0.
\ee
We can write now
\begin{align}
j^\epsilon S_0 w_j&=\sum_k \frac {e^{o_{\min \{j,k\}}(1)}
k^{\balpha-\epsilon} j^{\balpha-1+\epsilon}}
{k^{2\balpha}+2 \cos \theta k^\balpha j^\balpha+j^{2\balpha}}\\
\nonumber
&=\sum_{k \leq \ln j} \frac {e^{o_k(1)}
k^{\balpha-\epsilon} j^{\balpha-1+\epsilon}}
{k^{2\balpha}+2 \cos \theta k^\balpha j^\balpha+j^{2\balpha}}+
\sum_{k>\ln j} \frac {e^{o_j(1)}
k^{\balpha-\epsilon} j^{\balpha-1+\epsilon}}
{k^{2\balpha}+2 \cos \theta k^\balpha j^\balpha+j^{2\balpha}}.
\end{align}
Moreover,
\be
\sum_{k \leq \ln j} e^{o_k(1)}
\frac {k^{\balpha-\epsilon} j^{\balpha-1+\epsilon}}
{k^{2\balpha}+2 \cos \theta k^\balpha j^\balpha+j^{2\balpha}}=o_j(1),
\ee
provided $\epsilon<\balpha+1$
(for $\epsilon \geq \balpha+1$ the sum is not even $O_j(1)$), and
\be
\sum_{k>\ln j} \frac {e^{o_j(1)}
k^{\balpha-\epsilon} j^{\balpha-1+\epsilon}}
{k^{2\balpha}+2 \cos \theta k^\balpha j^\balpha+j^{2\balpha}}=e^{o_j(1)}
\int_0^\infty \frac {t^{\balpha-\epsilon} j^\balpha} {t^{2\balpha}+2 \cos
\theta t^\balpha j^\balpha+j^{2\balpha}} dt=e^{o_j(1)} S_\epsilon,
\ee
provided that $|\epsilon-1|<\balpha$
(if $|\epsilon-1| \geq \balpha$ the sum is not even $O_j(1)$).
We can now conclude, for $|\epsilon-1|<\balpha$,
\be
w_j j^\epsilon=e^{o_j(1)} \frac {S_\epsilon} {S_0}
\ee
and for $|\epsilon-1| \geq \balpha$,
\be
\lim_{j \to \infty} w_j j^\epsilon=\infty.
\ee
In particular, $DT(X)$ is a bounded operator in $u(\epsilon)$
if and only if $|\epsilon-1|<\balpha$, in which case the bound is uniform on
$X \in K$.  Moreover, if
$0<\epsilon<2$ then there exists $S_\epsilon S^{-1}_0<\hat
C_\epsilon<1$ and $N>0$ (independent of $X \in K$) such
that for $j>N$,
\be
w_j j^\epsilon<\hat C_\epsilon.
\ee
Let us now fix $N$ as above.
Let $v'_k=\min \{N^{-\epsilon},k^{-\epsilon}\}$, and $w'=DT(X) v'$.
\comm{
The stochasticity of $DT$ gives clearly
\be \label {improve}
\|w'\|_0 \leq \|v'\|_0.
\ee
Since
\be
|v'_{N+1}|=\frac {N^\epsilon} {(N+1)^\epsilon} \|v'\|_0<\|v'\|_0,
\ee
the positivity of $DT$ allows us to improve (\ref {improve}) to
\be
|w'_k|<\|v'\|_0, \quad k \geq 1,
\ee
where the inequality is uniform on $X \in K$ provided we take $k$ in a
bounded range.  In particular there exists $\tilde C_\epsilon<1$ with
\be
\sup_{k \leq N} w'_k \leq \tilde C_\epsilon N^{-\epsilon}
\ee
and indeed $\tilde C_\epsilon$ is independent of $X \in K$.}
By Lemma \ref {ineq}, we have
\be
\|w'\|_0<\|v'\|_0=N^{-\epsilon},
\ee
where the inequality is uniform on $X \in K$ (using for instance Lemma \ref
{lipderiv}), so there exists $\tilde C_\epsilon<1$ independent of $X \in K$
with
\be
\sup_{k \leq N} w'_k \leq \tilde C_\epsilon N^{-\epsilon}.
\ee

Let
\be
\|u\|_c=\sup_k \frac {|u_k|} {|v'_k|},
\ee
which is equivalent to the usual norm on
$u(\epsilon)$, since
\be
\|u\|_\epsilon \leq \|u\|_c \leq N^\epsilon \|u\|_\epsilon.
\ee
Clearly $\|DT(X) u\|_c \leq C_\epsilon \|u\|_c$ with
$C_\epsilon=\max\{\hat C_\epsilon, \tilde C_\epsilon\}$.
\end{pf}

\begin{rem}

Let $v_k=k^{-\epsilon}$ and $v^{(n)}=DT^n(P) v$.
A lower bound for the spectral radius of $DT(P)$
in $u(\epsilon)$ is given by
\be
\limsup_{n \to \infty} \|DT^n(P) v\|^{1/n} \geq \lim_{n \to \infty}
\lim_{k \to \infty} (k^\epsilon v^{(n)}_k)^{1/n}=\frac {S_\epsilon} {S_0}.
\ee
This achieves a minimum at $\epsilon=1$ and one
actually has $S_1 S_0^{-1}=\balpha-1$ (see Remark~\ref {r2}).
Notice that as $\theta \to \pi$
(which happens when $M \to \infty$ for the anharmonic oscillator),
$\balpha=1+\frac {\theta} {\pi} \to 2$ so
the contraction factor becomes weak.  This
should be compared to numerical estimates in \cite {V}, \S 7.3.

\end{rem}

\comm{
Let $X \in U(\balpha,\balpha)$ and for $\beta<\balpha$
let $U(\balpha,\beta,X)$ be the space of $Y$ with
\be
y_k=x_k+O(k^{\beta-\balpha}).
\ee
}

\begin{cor} \label {cor3}

Let $X \in U^0(\balpha,0)$.  If $0<\epsilon<\balpha+1$ then
$T(U(X,\epsilon))=U(T(X),\epsilon)$.

\end{cor}

\begin{pf}

Integrate the previous estimate.
\end{pf}

\begin{cor} \label {cor4}

If $0<\epsilon<2$ then $P$ is a global exponential attractor
in $U(P,\epsilon)$.

\end{cor}

\begin{pf}

Let $\underline X \leq P \leq \overline X \in U(P,\epsilon)$, and let
$K=\{\underline X \leq Y \leq \overline X\}$.  Then $\cup_{n=0}^\infty
T^n(K)$ is uniformly bounded in $U^0(\balpha,0)$ (Lemma \ref {precompact}),
and by Lemma \ref {strong}
there exists $C<1$ and a norm $\| \cdot \|_c$ in $u(\epsilon)$ such that if
$X \in K$ then $\|DT^n(X) v\|_c \leq C^n \|v\|_c$.  Integrating this
inequality we see that if $Y \in K$ then $\|T^n(y)-p\|_c \leq C^n
\|y-p\|_c$.
\end{pf}

Theorem \ref {actual result} follows from
Corollaries \ref {cor1}, \ref {cor3}, \ref {cor4} and Lemmas \ref {drif} and
\ref {lemma2}.

\begin{rem}

Let us remark that while the operator $T$ in $U(\balpha,0)$ has a line of
fixed points $\lambda P$, $\lambda>0$ (where $P$ is the fixed point in
$U^0(\balpha,0)$), this line is not a global attractor
in the $U(\balpha,0)$ metric.  Indeed it is easy to see that if
$1=n_1<n_2<...$ is a sequence that grows sufficiently fast and
\be
v_k=-1, n_{2j-1} \leq k <n_{2j}, \quad v_k=1, n_{2j} \leq k <n_{2j+1},
\ee
then letting $x=p+v$, $x^{(n)}=T^n(x)$ we have
\be
\inf_{\lambda>0} \|x^{(n)}-\lambda p\|_0=\|x^{(n)}-p\|_0=1,
\ee
for all $n \geq 0$, and we do not even have pointwise convergence:
\be
\liminf_{n \to \infty} x^{(n)}_k-p_k=-1,\quad \limsup_{n \to \infty}
x^{(n)}_k-p_k=1,
\ee
for all $k \geq 1$.

\end{rem}

\begin{rem}

A construction similar to the previous remark shows that $P$ is far from
being exponentially attractive in the $U^0(\balpha,0)$ metric: for any
decreasing sequence $a_1>a_2>...$ with $\lim_{k \to \infty} a_k=0$,
there exists $X \in U^0(\balpha,0)$ such that $\|T^n(x)-p\|_0>a_n$.

\end{rem}

{\bf Acknowledgements:}  I would like to thank Andr\'e Voros and
Jean-Christophe Yoccoz for several useful discussions, which originated
many of the arguments in this paper.

%%% Local Variables: 
%%% mode: latex
%%% TeX-master: t
%%% End: 

\end{document}